\documentclass[12pt]{amsart}
\usepackage[mathscr]{eucal}

\newtheorem{theorem}{Theorem}[section]
\newtheorem{lemma}[theorem]{Lemma}
\newtheorem{proposition}[theorem]{Proposition}
\newtheorem{corollary}[theorem]{Corollary}

\theoremstyle{definition}
\newtheorem{definition}[theorem]{Definition}

\newtheorem{example}[theorem]{Example}
\newtheorem{question}[theorem]{Question}

\theoremstyle{remark}
\newtheorem*{remark}{Remark}

\title{Abundant configurations in sumsets with one dense summand}
\author{John  T. Griesmer}
\date{\today}

\begin{document}

\begin{abstract}{We analyze sumsets $A+B$ $(=\{a+b:a\in A, b\in B\})$ where $A$ and $B$ are sets of integers, $A$ is infinite, and $B$ has positive upper Banach density.  For each $k,$ we show that $A+B$ contains at least the expected density of $k$-term arithmetic progressions, based on the density of $B,$ in contrast with an example of Bergelson, Host, Kra, and Rusza.  Furthermore, we show that when $A$ is infinite and $B$ has positive upper Banach density, $A+B$ must contain finite configurations not found in arbitrary sets of positive density, in contrast with results of Frantzikinakis, Lesigne, and Wierdl on sets of $k$-recurrence.}
\end{abstract}

\maketitle

\section{Introduction}

A recurring theme in additive combinatorics is that sumsets $A+B:=\{a+b:a\in~A, b\in B\}$ tend to be more structured than their summands $A$ and $B.$  Bourgain (\cite{Bourgain}) demonstrated a striking example of this phenomenon: if $A,B\subset [1,\dots, N]$ with $|A|=\alpha N, |B|=\beta N,$ then for sufficiently large $N$ (depending on $\alpha, \beta$) $A+B$ contains an arithmetic progression of length approximately $e^{\alpha \beta (\log N)^{1/3}}.$  This result was strengthened and generalized in \cite{BG}, \cite{HL}, \cite{San}, and \cite{CS}.  

\medskip

We will be concerned with infinite sets of integers; Section \ref{integers} is a glossary of relevant terminology.  R.~Jin showed that when $A,B\subset \mathbb Z$ both have positive upper Banach density, $A+B$ must be piecewise syndetic (\cite{Jin1}).  This result was strengthened and generalized in \cite{JK},\cite{BFW},\cite{BBF}, \cite{Beiglbock}, and \cite{Griesmer}.  A consequence of the main theorem of \cite{BFW} is that the density of $k$-term arithmetic progressions (with certain common differences) in $A+B$ is close to the maximum of the densities of $A$ and $B.$  Furthermore, $A+B$ will contain certain kinds of configurations which are not necessarily present in sets of positive density.   These conclusions holds even when $A$ is assumed to satisfy a weaker hypothesis than positive density, as shown in \cite{Griesmer}.

\medskip

As mentioned above, sets $B\subset \mathbb Z$ of positive density need not be as structured as sumsets $A+B.$  In particular, $B$ need not be piecewise syndetic.  Nevertheless, classical and recent results indicate that sets of positive density contain highly structured subsets.  In particular, Szemer\'edi's Theorem on arithmetic progressions guarantees that every set of positive upper Banach density contains arithmetic progressions of every finite length.

\begin{theorem}[Szemer{\'e}di's Theorem, \cite{Sz}] \label{SzThm}  Let $E\subset \mathbb Z$ with $d^*(E)>0.$  For all $k\in \mathbb N,$ there exists $n\in \mathbb N$ such that $E\cap E-n \cap E-2n \cap \dots \cap E-(k-1)n\neq \emptyset.$
\end{theorem}
The methods of \cite{BHMP} show how Furstenberg's ergodic-theoretic proof (\cite{F77}) of Theorem \ref{SzThm} gives a qualitatively stronger conclusion:
\begin{theorem}\label{MRThm} For all $k, \delta>0,$ there exists $C=C(k,\delta)>0$ such that
$$
\{n:d^*(E\cap E-n \cap E-2n \cap \dots \cap E-(k-1)n)> C(k,\delta)-\varepsilon\}
$$
is syndetic whenever $d^*(E)\geq \delta, \varepsilon>0.$
\end{theorem}

Naturally one wonders: what are the optimal values of $C(k,\delta)$ in Theorem \ref{MRThm}?  We will a give standard example in Section \ref{taxonomy} showing that $C(k,\delta)\leq \delta^k$ for each $k\geq 1.$  The following result from \cite{BHK} shows that $C(k,\delta)=\delta^k$ for $k=2, 3,$ and $4,$ while $C(k,\delta)$ is much smaller than $\delta^k$ for $k>5.$

\begin{theorem}\label{BHKThm} (1) Let $E\subset \mathbb Z$ with $d^*(E)>0.$  Then for all $\varepsilon>0,$ the sets
$$
\{n: d^*(E\cap E-n \cap E-2n)> d^*(E)^3 -\varepsilon\},
$$
$$
\{n: d^*(E\cap E-n \cap E-2n\cap E-3n)> d^*(E)^4 -\varepsilon\},
$$
are syndetic.

(2)  For all $l>0,$ there exists $E\subset \mathbb Z$ with $d^*(E)>0$ such that
$$
d^*(E\cap E-n \cap E-2n \cap E-3n \cap E-4n)< d^*(E)^l/2
$$
for all $n\neq 0.$
\end{theorem}

The set $R_k(E):=\{n:d^*(E\cap E-n \cap E-2n \cap \cdots \cap E-(k-1)n)>0\}$ has received considerable study.  For instance, it is known that for all $k, R_k(E)$ contains perfect squares.  S\'ark{\H{o}}zy (\cite{Sar}), and Furstenberg (\cite{F77}) independently proved this in the case $k=2.$  For $k>2$ this is due to Bergelson and Leibman (\cite{BL}).  The example $E=4\mathbb Z$ shows that $R_2(E)$ need not contain numbers of the form $n^2+1,$ as every square is congruent to $0$ or $1 \text{ mod }4$.  More generally, a set $S\subset \mathbb Z$ is called \textit{$k-1$-intersective} if for all $E\subset \mathbb Z$ with $d^*(E)>0,$ there exists $n\in S$ such that $d^*(E\cap E-n\cap E-2n \cap \dots \cap E-(k-1)n)>0.$  Furstenberg (\cite{Fbook}) showed that there is a $1$-intersective set $S$ which is not $2$-intersective, while \cite{FLW} demonstrates the existence of $k-1$-intersective sets $S$ that are not $k$-intersective for all $k\geq 2.$  The main theorem from \cite{FLW} is as follows.

\begin{theorem}\label{FLWThm} For irrational $\alpha\in \mathbb R$ and $k\geq 2,$ the set $S_{k,\alpha}:=\{n:n^k\alpha\mod 1\in (\frac{1}{4},\frac{3}{4})\}$ is $k-1$-intersective but not $k$-intersective.
\end{theorem}
In particular, for each irrational $\alpha$ and $k\geq 2$ there is a set $E\subset \mathbb Z$ with $d^*(E)>0$ and $E\cap E-n\cap E-2n\cap \dots \cap E-(k-1)n=\emptyset$ for all $n\in S_{k,\alpha}.$

\medskip

This article will study the kinds of configurations found in sets of the form $A+B,$ where $A\subset \mathbb Z$ is infinite and $B\subset \mathbb Z$ has positive upper Banach density.  Here is a special case of our main result, which contrasts sharply with Part 2 of Theorem \ref{BHKThm} and Theorem \ref{FLWThm}.

\begin{theorem}\label{special}  Let $A\subset \mathbb Z$ be infinite and let $d^*(B)>0.$  Then for all $k\in \mathbb N$ and all $\varepsilon>0,$ there exist $n\in \mathbb N$ with
\begin{align}\label{expected}
d^*(A+B\cap A+B-n\cap A+B-2n\cap \dots \cap A+B-(k-1)n)\geq d^*(B)^k-\varepsilon.
\end{align}
For all irrational $\alpha,$ there exist $n\in S_{\alpha,j}:=\{m: m^j\alpha\mod 1\in (1/4,3/4)\}$ satisfying (\ref{expected}), and the set of such $n$ is syndetic.
\end{theorem}

In Section \ref{proofs} we summarize some background and results from ergodic theory.  We will deduce Theorem \ref{special} from Theorem \ref{main}, an analogous result about measure preserving systems.  Theorem \ref{main} will be a fairly straightforward consequence of the main result of \cite{BHK}, together with a classical result of Leon Green (Lemma \ref{lebesguespectrum}) concerning the spectral type of nilsystems.

\subsection{Stronger hypotheses on $A$}

Supposing both $A,B\subset \mathbb Z$ have positive upper Banach density, the conclusion of Theorem \ref{main} can be strengthened considerably.  The following is implicit in the proof of \cite{BFW}, Theorem I.  An explicit proof appears in \cite{Griesmer}.

\begin{theorem}\label{max} Let $A,B,\subset \mathbb Z.$  Then for all $k\in \mathbb N, \varepsilon>0,$ there exists $n\in \mathbb N$ such that
$$
d^*(A+B\cap A+B-n\cap A+B-2n \cap \dots \cap A+B-(k-1)n)\geq \max\{d^*(A),d^*(B)\}-\varepsilon.
$$
In fact, the set of such $n$ is a Bohr set.
\end{theorem}

Our second result generalizes Theorem \ref{max} to cases where $d^*(A)=0, d^*(B)>0,$ including some not covered in \cite{Griesmer}.

\begin{definition}  Let $\mathcal S:=(S_j)_{j\in \mathbb N}$ be a sequence of finite subsets of $\mathbb Z.$  We say that $\mathcal S$ is \textit{nearly equidistributed} if for all all but countably many $\theta\in (0,2\pi),$
$$
\lim_{j\to \infty} \frac{1}{|S_j|}\sum_{n\in S_j} e^{in\theta}=0.
$$
If the above holds for all $\theta\in (0,2\pi),$ we say that $\mathcal S$ is \textit{equidistributed}.
\end{definition}

\begin{example}  The following $S_j$ form nearly equidistributed, but not equidistributed sequences.
\begin{enumerate}
\item[(i)] $S_j=\{1,4,9,16,\dots, j^2\}.$
\item[(ii)] $S_j=\{p\leq j: p \text{ is prime}\}.$
\end{enumerate}
The following $S_j$ form equidistributed sequences.
\begin{enumerate}
\item[(iii)] $S_j=\{2^j,2^j+1,2^j+2,\dots, 2^j+j-1\}.$
\item[(iv)] $S_j=\{\lfloor 1^{5/2} \rfloor,\lfloor 2^{5/2} \rfloor,\lfloor 3^{5/2} \rfloor,\dots, \lfloor j^{5/2} \rfloor\}.$
\item[(v)] $S_j=\{\lfloor \sqrt{2}j^5-\pi j^3\rfloor: 1\leq m\leq n\}.$
\end{enumerate}
Here $\lfloor x\rfloor$ denotes the greatest integer less than or equal to $x.$
\end{example}

Example (i) is a special case of Weyl's theorem on equidistribution, (ii) is due to Vinogradov,  (iii) is classical, (iv) is found in \cite{Wierdl}, and (v) is found in \cite{BKQW}, which has an extensive bibliography and is otherwise an excellent introduction and reference.

If $\mathcal S=(S_n)_{n\in \mathbb N}$ is a sequence of sets of integers and $A\subset \mathbb Z,$ the \textit{upper density of $A$ relative to $\mathcal S$} is
$$
d_{\mathcal S}:=\limsup_{n\to \infty} \frac{|A\cap S_n|}{|S_n|}.
$$

\begin{theorem}\label{near}  Suppose that $\mathcal S$ is a nearly equidistributed sequence, $A,B\subset \mathbb Z$ with $d_{\mathcal S}(A)>0, d^*(B)>0.$  Then for all $\varepsilon>0$ and all $k\in \mathbb N,$ there exists $n\in \mathbb N$ such that
$$
d^*(A+B\cap A+B-n \cap A+B-2n \cap \dots \cap A+B-(k-1)n)\geq d^*(B)-\varepsilon.
$$
In fact, there is a Bohr set of such $n.$
\end{theorem}

In the next section we consider these issues from an ergodic-theoretic perspective.

\subsection{Acknowledgements.}  The author is pleased to thank Izabella {\L}aba for financial support, Izabella {\L}aba and Malabika Pramanik for helpful discussions, and Vitaly Bergelson for clarifying the relationship between this work and the work in \cite{BHMP}.

\bigskip

\section{Main results}\label{proofs}

Our results about sets of integers will be deduced from corresponding results in ergodic theory.  First we outline some notation, terminology, and background, then we discuss some special classes of dynamical systems and some special classes of sets of integers.

\subsection{Measure preserving systems}

A \textit{measure preserving system}, or M.~P.~S., is a quadruple $(X,\mathscr X,\mu,T)$, where $(X,\mathscr X,\mu)$ is a measure space with probability measure $\mu,$ and $T:X\to X$ is a map preserving $\mathscr X$ and $\mu$ in the sense that $T^{-1} \mathscr X\subset \mathscr X$ and $\mu(T^{-1}D)=\mu(D)$ for all $D\in \mathscr X.$  We will assume $T$ is invertible.
\medskip

\textit{Notation}:  The symbol $\mathbf X$ will always be an abbreviation of the sequence of symbols $(X,\mathscr X,\mu,T),$ similarly $\mathbf Y$  will always abbreviate $(Y,\mathscr Y, \nu,S).$

\medskip

If $\mathbf X =(X,\mathscr X,\mu, T)$ is a M.~P.~S., a \textit{factor} of $\mathbf X$ is a M.~P.~S. $\mathbf Y=(Y,\mathscr Y,\nu, S)$ together with a map $\pi:X\to Y$ satisfying (i) $\pi^{-1}\mathscr Y\subset \mathscr X,$ (ii)  $\mu(\pi^{-1}D)=\nu(D)$ for all $D\in \mathscr Y,$ and (iii) $S\pi(x)=\pi(Tx)$ for $\mu$-almost every $x\in X.$  We write $\pi: \mathbf X\to \mathbf Y$ to specify that $\mathbf Y$ is a factor of $\mathbf X$ with factor map $\pi.$  If $\pi:X\to Y$ is one-to-one on a set of full measure, then we say that $\pi$ is an isomorphism and $\mathbf X$ and $\mathbf Y$ are isomorphic.

\medskip

Every factor $\mathbf Y$ of $\mathbf X$ determines a $T$-invariant $\sigma$-algebra $\pi^{-1}\mathscr Y\subset \mathscr X,$ and conversely, associated to every $T$-invariant $\sigma$-algebra $\mathscr B\subset \mathscr X,$ there is a factor $\mathbf Y$ such that $\pi^{-1}\mathscr Y=\mathscr B,$ up to $\mu$-measure $0.$  Given a factor $\pi:\mathbf X\to \mathbf Y$ we may abuse notation and write $\mathscr Y$ for $\pi^{-1}\mathscr Y.$  We may then identify $L^2(\nu)$ with the space of $\pi^{-1}\mathscr Y$-measurable square-integrable functions, and consider the conditional expectation map $f\mapsto \mathbb E(f|\pi^{-1}\mathscr Y),$ or $f\mapsto \mathbb E(f|\mathscr Y)$ with our abuse of notation.  When $f\in L^2(\mu), \mathbb E(f|\mathscr Y)$ agrees with the orthogonal projection of $f$ on the closed subspace of $\mathscr Y$-measurable functions in $L^2(\mu).$

\medskip

Supposing $\pi_i:\mathbf X\to \mathbf Y_i$ is a directed collection of factors of $\mathbf X=(X,\mathscr X,\mu, T),$ we say that $\mathbf X$ is an \textit{inverse limit} of the $\mathbf Y_i$ if the $\sigma$-algebras $\pi_i^{-1} \mathscr Y_i$ generate $\mathscr X$ (up to $\mu$-measure $0$).

\medskip

A system $\mathbf X=(X,\mathscr X,\mu,T)$ is called \textit{ergodic} if $\mu(T^{-1}D)=\mu(D)$ implies $\mu(D)=0$ or $\mu(D)=1.$  Equivalently, $\mathbf X$ is ergodic if the only functions $f\in L^2(\mu)$ satisfying $f\circ T=f$ are constant $\mu$-almost everywhere.





\subsection{Taxonomy of measure preserving systems}\label{taxonomy}

Here we describe the special classes of systems we consider in the sequel.

\medskip

\noindent \textbf{Kronecker systems.}  If $Z$ is a compact abelian group and $\alpha\in Z$ is such that $\{n\alpha: n\in \mathbb Z\}$ is dense in $Z,$ consider the map $R_\alpha:Z\to Z$ given by $R_\alpha(z)=z+\alpha.$  Then $R_\alpha$ preserves the Haar measure $m$ of $Z,$ and $(Z,\mathscr Z,m,R_\alpha)$ is an ergodic M.~P.~S.  Such systems are called \textit{Kronecker systems}.  Kronecker systems have \textit{discrete spectrum}, meaning $L^2(m)$ is spanned by the eigenvectors of $T$ - those functions satisfying $f\circ T=\lambda f$ for some $\lambda \in \mathbb C.$  The Halmos-von Neumann theorem (\cite{Hal} or \cite{Glasner})
says that every ergodic M.~P.~S.~having discrete spectrum is isomorphic to a Kronecker system.

\medskip

Every M.~P.~S.~ $\mathbf X$ has a maximal factor with discrete spectrum, namely, the factor corresponding to the $\sigma$-algebra generated by the eigenvectors of $T.$  If $\mathbf X$ is ergodic, this factor is called the \textit{Kronecker factor}, and is isomorphic to a Kronecker system, by the Halmos-von Neumann theorem.

\medskip

 \noindent \textbf{Nilsystems.}  Fix $k,$ and let $G$ be a $k$-step nilpotent Lie group, and $\Gamma\subset G$ a discrete subgroup such that $G/\Gamma$ is compact.  Then $G/\Gamma$ has a  Borel probability measure $\mu$ invariant under the action of $G,$ that is $\mu(g\cdot D\Gamma)=\mu(D\Gamma)$ for every Borel set $D\subset G.$  If $g\in G,$ we consider the M.~P.~S.~$(X,\mathscr X,\mu,T_g),$ where $X=G/\Gamma, \mathscr X$ is the Borel $\sigma$-algebra of $G/\Gamma,$ and $T_gx=g\cdot x.$  Such a M~.P.~S.~ is called a \textit{$k$-step nilsystem}, or simply a \textit{nilsystem}.

\medskip

The correlation sequences of nilsystems are particularly nice;  the following description of such sequences is essentially due to Leon Green.
 \begin{lemma}\label{lebesguespectrum}
 If $(X,\mathscr X,\mu,T)$ is a nilsystem with Kronecker factor $(Z,\mathscr Z,m,R_\alpha),$ and $f\in L^2(\mu)$ with $\mathbb E(f|\mathscr Z)=0,$ then for all $g\in L^2(\mu),$
$$
\lim_{n\to \infty} \int f\circ T^n \cdot g\, d\mu=0.
$$
\end{lemma}
This is a slightly different formulation of what is proved in \cite{AGH}, Chapter 5, and technically Green's result does not imply the above statement.  A result of Parry (\cite{Parry}) does imply the above statement; we outline the deduction of Lemma \ref{lebesguespectrum} from Parry's result in Section \ref{appendix}.

\medskip

A special class of nilsystems is formed by taking repeated group extensions of torus rotations.  If $Z=\mathbb T^d$ is the $d$-dimensional torus and $\alpha\in Z$ is such that $\{n\alpha:n\in \mathbb N\}$ is dense in $Z,$ then $(Z,\mathscr Z,m,R_\alpha)$ is an ergodic M.~P.~S.  If $\psi: Z\to Z$ is an affine map (a homomorphism plus a constant) we may form a new M.~P.~S.
$(Z\times Z, \mathscr Z\otimes \mathscr Z, m\times m, T_{\alpha, \psi})$ where $T_{\alpha,\psi}(x,y)=(x+\alpha,y+\psi(x)).$  For instance, when $Z=\mathbb T$ and $\psi(x)=2x+\alpha,$ we have $T^n(0,0)=(n\alpha, n^2\alpha).$

If $p$ is a polynomial of degree $d,$ this construction may be iterated, as shown in \cite{Fbook}, to obtain an ergodic system on $Z^d$ such that $T^n(0,\dots, 0)$ has $p(n)\alpha$ as its last coordinate.  It is known that the resulting system $(Z^d, \mathscr Z^{\otimes d}, m^{d},T)$ is actually a $d$-step nilsystem, see, for instance \cite{BLpoly}.

\medskip

\noindent \textbf{Weakly mixing systems.}   A system $\mathbf X$ is said to be weakly mixing if the only eigenfunctions of $T$ are constant almost everywhere.  Furstenberg (\cite{F77})  showed that $\mathbf X$ is weakly mixing if and only if for all $k\in \mathbb N$ and all $f\in L^2(\mu),$
$$
\lim_{N-M\to \infty} \frac{1}{N-M}\sum_{n=M}^{N-1} \bigl| \int f\cdot f\circ T^n\cdot f\circ T^{2n} \cdot \cdots \cdot f\circ T^{(k-1)n}\, d\mu - \Bigl( \int f\, d\mu \Bigr)^k\bigr|=0,
$$
or equivalently that for all $\varepsilon>0,$
$$
d^*\Bigl\{n: \bigl| \int f\cdot f\circ T^n\cdot f\circ T^{2n} \cdot \cdots \cdot f\circ T^{(k-1)n}\, d\mu - \Bigl( \int f\, d\mu \Bigr)^k\bigr|>\varepsilon \Bigr\}=0.
$$
Bergelson (\cite{PET}) extended this characterization to show that $\mathbf X$ is weakly mixing if and only if for every $k$ and every $k-1$-tuple of polynomials $p_1,\dots, p_{k-1},$
$$
d^*\Bigl\{n: \bigl| \int f\cdot f\circ T^{p_1(n)}\cdot f\circ T^{p_2(n)} \cdot \cdots \cdot f\circ T^{p_{k-1}(n)}\, d\mu - \Bigl( \int f\, d\mu \Bigr)^k\bigr|>\varepsilon \Bigr\}=0.
$$

\medskip

\noindent \textbf{Rigid systems.}  A system $\mathbf X$ is called \textit{rigid} if there is an increasing sequence $\{n_j\}_{j\in \mathbb N}$ satisfying $\lim_{n\to \infty} f\circ T^n=f$ weakly for all $f\in L^2(\mu).$  There are weakly mixing, rigid systems $\mathbf X,$ and with the appropriate topology on the set $G$ of invertible Lebesgue-measure preserving transformations of $[0,1],$ the set of weakly mixing, rigid systems is a dense $G_\delta$ in $G.$  See \cite{Walterspaper} for details.  The existence of such a system leads to the following examples, showing that our main theorem is sharp in one sense.

\begin{example}\label{optimal}  Let $\mathbf X$ be a rigid, weakly mixing M.~P.~S., $D\in \mathscr X,$ and $p_i:\mathbb Z\to \mathbb Z,i=1,\dots, k$ polynomials.  For all $\varepsilon>0,$ there exists $A\subset \mathbb Z$ such that $D_A:=\bigcup_{a\in A} T^{a}D$ satisfies $\mu(D_A)<\mu(D)+\varepsilon/2.$  Since $\mathbf X$ is weakly mixing, we have
$$
d^*\bigl\{n: \mu\bigl(D_A \cap \bigcap_{i=1}^{k-1} T^{-p_i(n)} D_A\bigr)>\mu(D_A)^k+\varepsilon\}=0.
$$
In particular, the above set of integers is not syndetic.
\end{example}

\medskip

\subsection{Sets of integers.}\label{integers}  Here we summarize the various properties of sets of integers recently studied; the terminology is mainly to describe the sets of differences $n$ when the intersection $d^*(E\cap E-n\cap \dots \cap E-(k-1)n)$ is large, if $E\subset \mathbb Z$ is a set with some special structure.

\medskip

For $S\subset \mathbb Z,$ the \textit{upper Banach density} is the number $d^*(S):=\lim_{M\to \infty} \sup_{N\in \mathbb Z} \frac{|A\cap [N,N+M]|}{N+M-1}.$

\medskip

A set $S\subset \mathbb Z$ is called \textit{syndetic}, if one of the following equivalent conditions holds: (i) there exists a finite $F\subset \mathbb Z$ such that $F+S=\mathbb Z,$ (ii) $S\cap R\neq \emptyset$ whenever $R$ is a set containing arbitrarily long finite intervals.

\medskip

$R$ is called \textit{thick} if $R$ contains intervals of every finite length, and $S$ is called \textit{piecewise syndetic} if $S=S'\cap R$ where $S'$ is syndetic and $R$ is thick.

\medskip

Examples of syndetic sets are given by certain dynamical systems.  A topological system $(X,T)$ is called \textit{minimal} if $X$ contains no closed $T$-invariant subset.  For such $(X,T),$ the set $\{n:T^nx\in U\}$ is syndetic whenever $U\subset X$ is open and $x\in X.$  When viewed as topological systems, Kronecker systems and nilsystems are minimal, and the sets of entry times arising therefrom have distinguished properties.  Following \cite{BFW}, $S\subset \mathbb Z$ is called a \textit{Bohr set} if there is a Kronecker system $(Z,\mathscr Z,m,R_\alpha)$ and a set $U\subset Z$ such that $\{n: n\alpha\in U\}\subset S.$  Following \cite{BFW}, $S$ is called \textit{piecewise Bohr} if $S=S'\cap R,$ where $S'$ is Bohr and $R$ is thick.

\medskip

Following \cite{HKnb}, $S$ is called \textit{Nil-Bohr} if there is a nilsystem $\mathbf Z=(Z,\mathscr Z, \mu_Z,T)$ and an open $U\subset Z$ and $x\in X$ with $\{n:T^nx\in U\}\subset S,$ and $S$ is called \textit{piecewise Nil-Bohr} if $S=S'\cap R$ where $S'$ is Nil-Bohr and $R$ is thick.  If $\mathbf Z$ is a $k$-step nilsystem, then ``Nil-Bohr" may be specialized to ``$\text{Nil}_k$-Bohr."

\medskip

A special class of $\text{Nil}_k$-Bohr sets is given by polynomial orbits of Kronecker systems: if $\mathbf Z$ is a Kronecker system and $p:\mathbb Z\to \mathbb Z$ is a polynomial with $p(0)=0,$ then for all neighborhoods $U$ of $0,$ the set $R_{p,U}:=\{n:p(n)\alpha\in U\}$ is a $\text{Nil}_k$-Bohr set  (see \cite{BLpoly}).  Furthermore, the intersection of finitely many such $R_{p,U}$ is again an $\text{Nil}_k$-Bohr set.  This leads to the following multiple recurrence result for Kronecker systems.

\begin{theorem}\label{grouprecurrence}  Let $(Z,\mathscr Z,m,R_\alpha)$ be a Kronecker system, and let $p_1,\dots, p_{k-1}:\mathbb Z\to \mathbb Z$ be polynomials of degree at most $d.$  For all $D\subset Z$ and all $\varepsilon>0,$
$$
\{n: \mu(D\cap D-p_1(n)\alpha \cap D-p_2(n)\alpha \cap \dots \cap D-p_{k-1}(n)\alpha)> \mu(D)-\varepsilon\}
$$
is $\text{Nil}_d$-Bohr.
\end{theorem}

\medskip

Note that a (piecewise) syndetic set $S$ remains such after removing a set of upper Banach density zero: $S\setminus A$ is (piecewise) syndetic if $d^*(A)=0.$ We say that $S$ is \textit{almost} Bohr (or \textit{almost} $\text{Nil}_k$-Bohr) if $S= S'\setminus A,$ where $S'$ is Bohr (or $\text{Nil}_k$-Bohr), and $d^*(A)=0.$

\subsection{The correspondence principle}  Furstenberg (\cite{F77}) established a general principle relating finite configurations in sets of positive density to measure preserving systems.  We use a version of this principle from \cite{BHK}, modified slightly to deal with sumsets.

\begin{proposition}[cf.~ \cite{BHK}, Proposition 3.1] \label{correspondence} Let $B\subset \mathbb Z.$  Then there is an ergodic M.~P.~S. $\mathbf X=(X,\mathscr X,\mu,T)$ with $D\in \mathscr X$ satisfying $\mu(D)=d^*(B)$ and
\begin{align*}
\mu(T^{-n_1}D_A \cap T^{-n_2}D_A \cap \cdots \cap T^{-n_k}D_A) \leq d^*( A+B-n_1\cap A+B-n_2\cap \dots \cap A+B-n_k)
\end{align*}
where $D_A=\bigcup_{a\in A} T^{a} D,$ for all finite $A\subset \mathbb Z$ and all $n_1,\dots, n_k\in \mathbb Z.$
\end{proposition}

\textit{Proof.}  Following \cite{BHK}, let $X\subset \{0,1\}^\mathbb Z$ be the orbit closure of $1_B$ under the shift $T,$ where $(Tx)(n)=x(n+1).$  Let $D=\{x\in X: x(0)=1\},$ and repeat the steps of \cite{BHK}, Proposition 3.1 with $D_A= \bigcup_{a\in A}T^a D$ replacing $D,$ noting that the indicator function $1_{D_A}$ is continuous. \hfill $\square$

\subsection{Sets of $k$-recurrence and $k$-intersective sets}

As presented in \cite{F77}, \cite{FLW}, a set $S\subset \mathbb Z$ is called a \textit{set of $k-1$-recurrence} if for every M.~P.~S.~$\mathbf X,$ and every $D\in \mathscr X$ with $\mu(D)>0,$ there exists $n\in S$ such that
$$
\mu(D\cap T^{-n}D\cap T^{-2n}D\cap \dots \cap T^{-(k-1)}D)>0.
$$
Similarly, $S$ is called $k-1$-intersective if for all $B\subset \mathbb Z$ with $d^*(B)>0,$ there exists $n\in S$ with
$$
d^*(B\cap B-n\cap B-2n \cap \dots \cap B-(k-1)n)>0.
$$
It is shown in \cite{FKO} that for all $k,$ the $k-1$-intersective sets are exactly the sets of $k-1$-recurrence, so in the sequel we speak only of sets of $k-1$-recurrence.

\medskip

We may consider more general forms of recurrence.  Let $p_1,\dots, p_{k-1}:\mathbb Z\to \mathbb Z$ be functions, and write $\vec{p}:=(p_1,\dots, p_{k-1}).$  Call a set $S\subset \mathbb Z$ a \textit{set of $\vec{p}$-recurrence} if for every M.~P.~S.~$\mathbf X$ and every $D\in \mathscr X,$ there exists $n\in S$ such that $\mu(D\cap T^{-p_1(n)}D\cap T^{-p_2(n)}D\cap \dots \cap T^{-p_{k-1}(n)}D)>0,$ and call $\vec{p}$ recurrent if $\mathbb N$ is a set of $\vec{p}$-recurrence.

\medskip

Some variations on this definition will be useful.  Call $S$ a \textit{frequent set of $\vec{p}$-recurrence} if for every M.~P.~S. $\mathbf X$ and all $D\in \mathscr X$ with $\mu(D)>0,$ the set of $n$ satisfying $\mu(D\cap T^{-p_1(n)}D\cap T^{-p_2(n)}D\cap \dots \cap T^{-p_{k-1}(n)}D)>0$ has positive upper Banach density.

\medskip

Finally, call $S$ a \textit{set of $\vec{p}$-recurrence for Kronecker systems} if $S$ satisfies the definition of ``set of $\vec{p}$-recurrence" with ``M~.P.~S.~" replaced by ``Kronecker system."

\medskip

\begin{theorem}\label{main}  Let $p_i:\mathbb Z\to \mathbb Z$ be polynomials, $i=1,\dots, k-1.$  Suppose that $S$ is a set of frequent $\vec{p}$-recurrence for Kronecker systems.  Then for all ergodic M.~P.~S.~s $\mathbf X,$ all $D\in \mathscr X,$ all infinite $A\subset \mathbb Z,$ and all $\varepsilon>0,$
$$
R_\varepsilon:=\Bigl\{n: \mu(D_A\cap T^{-p_1(n)}D_A\cap T^{-p_2(n)}D_A\cap \dots \cap T^{-p_k(n)} D_A)>\mu(D_A)^k-\varepsilon \Bigr\}
$$
is almost $\text{Nil}_r$-Bohr where $r=\deg(\vec{p}),$ where $D_A=\bigcup_{a\in A}T^aD,$ furthermore, $d^*(S\cap R_\varepsilon)>0.$  If the $p_i$ are linear, then $R_\varepsilon$ is almost Bohr.
\end{theorem}

\medskip

It will be more convenient to consider functions rather than sets, so we reformulate Theorem \ref{main} as follows.

\medskip

\begin{theorem}\label{main'} Let $p_i:\mathbb Z\to \mathbb Z$ be polynomials, $i=1,\dots, k-1.$  Suppose that $S$ is a frequent set of $\vec{p}$-recurrence for Kronecker systems.  Then for every ergodic M.~P.~S.~$\mathbf X,$  $f:X\to [0,1],$ infinite $A\subset \mathbb Z,$  and all $\varepsilon>0,$ there exists a finite $A'\subset A$ such that
$$
R_\varepsilon:=\Bigl\{n:  \int f_{A'} \cdot f_{A'}\circ T^{p_1(n)} \cdot f_{A'}\circ T^{p_2(n)} \cdot \cdots \cdot f_{A'}\circ T^{p_{k-1}(n)}\, d\mu > \Bigl(\int f\, d\mu\Bigr)^k -\varepsilon \Bigr\}
$$
is almost $\text{Nil}_r$-Bohr, where $r=\deg(\vec{p}),$ where $f_{A'}=\frac{1}{|A'|}\sum_{a\in A'} f\circ T^{-a}.$  If the $p_i$ are linear, then $R$ is almost Bohr.  Furthermore, $d^*(S\cap R_\varepsilon)>0.$
\end{theorem}

Theorem \ref{main} follows from Theorem \ref{main'} by taking $f=1_D,$ and noting that $1_{D_{A'}}\geq f_{A'}$ whenever $A'\subset A.$

From Theorem \ref{main} and Proposition \ref{correspondence} we deduce the following combinatorial corollary.

\begin{corollary} Let $p_i:\mathbb Z\to \mathbb Z$ be polynomials, $i=1,\dots, k-1.$  Suppose that $S$ is a frequent set of $\vec{p}$-recurrence for Kronecker systems.  Then for every set $B\subset \mathbb Z$ with $d^*(B)>0,$
$$
R_\varepsilon:=\bigl\{n: d^*(A+B\cap A+B-p_1(n) \cap A+B-p_2(n) \cap \dots \cap A+B-p_{k-1}(n)) > d^*(B)^k-\varepsilon\bigr\}
$$
is almost $\text{Nil}_r$-Bohr, where $r=\deg(\vec{p}).$  If the $p_i$ are linear, then $R$ is almost Bohr.  Furthermore, $d^*(S\cap R)>0.$
\end{corollary}
We can deduce Theorem \ref{special} by noting that for each irrational $\alpha, S_{\alpha,k}\cap R$ is syndetic whenever $R$ is a Bohr set; this follows from Weyl's theorem on equidistribution, or from the topological recurrence methods in \cite{Fbook}.  Removing a set of upper Banach density $0$ will not affect syndeticity, so the conclusion remains even when $R$ is almost Bohr.

\begin{remark}  The systems of polynomials $(p_1,\dots,p_{k-1})$ for which $\mathbb N$ is a set of recurrence are characterized in \cite{BLL}
\end{remark}

\subsection{Questions}

\subsubsection{Optimality}  Can Theorem \ref{main} be improved?  Although Example \ref{optimal} shows that $$
R':=\{n:d^*(A+B\cap A+B-n\cap \dots \cap A+B-(k-1)n)\geq C\}
$$
is not necessarily syndetic for $C> d^*(B)^k,$ our methods do not resolve the following.

\begin{question}   If $A,B\subset \mathbb Z,$ with $A$ infinite, $d^*(B)>0,$  and $\varepsilon>0,$ does there exist $n\in \mathbb N$ with
$$
d^*(A+B\cap A+B-n\cap A+B-2n \cap \dots \cap A+B-(k-1)n)>d^*(B)-\varepsilon\text{?}
$$
If not, can $d^*(B)$ be replaced by some quantity larger than $d^*(B)^k$ to obtain an affirmative answer?

Similarly, if $\mathbf X$ is an ergodic M.~P.~S.~ and $D\subset X$ with $\mu(D)>0$ and $A\subset \mathbb Z$ is infinite, does there exist $n\in \mathbb N,$ with
$$
\mu(D_A\cap T^{-n}D_A\cap T^{-2n}D_A\cap \dots \cap T^{-(k-1)n}D_A)\geq \mu(D)-\varepsilon\text{?}
$$
\end{question}

\subsubsection{Generalizations to $\mathbb Z^d$}  For $S\subset \mathbb Z^d,$ one can define the upper Banach density $d^*(S):=\lim_{M\to \infty} \sup_{N\in \mathbb Z^d} \frac{|S\cap (N+I_M)|}{|I_M|},$ where $I_M$ is the cube $[0,M]^d.$  Many of the results about sumsets in $\mathbb Z$ generalize naturally to sumsets in $\mathbb Z^d,$ as in \cite{JK} and \cite{BBF}.  Similarly, Szemer\'edi's theorem generalizes to this setting, as shown \cite{FuKa}.  Our main result, however, breaks down under the naive generalization to $\mathbb Z^2:$ for all $l\in \mathbb N$ there exists $B\subset \mathbb Z^2$ with $d^*(B)>0$ and $A\subset \mathbb Z^2$ infinite with 
$$
d^*(A+B\cap A+B-(n,m)\cap A+B-2(n,m) \cap A+B-3(n,m) \cap A+B-4(n,m))\leq \frac{1}{2}d^*(B)^l
$$ 
for all $(n,m)\in \mathbb Z^2$ with $m\neq 0$:   take the example $B_0$ from \ref{BHKThm} and let $B=\mathbb Z\times B_0, A=\mathbb Z\times \{0\}.$

In this example we still have
$$
d^*(A+B-(n,0)\cap A+B-2(n,0)\cap \dots \cap A+B-4(n,0))= d^*(B)
$$
for all $n\in \mathbb Z,$ which raises the following question.

\begin{question}  Are there infinite sets $A,B\subset \mathbb Z^2$ with $d^*(B)>0$ such that for all $(n,m)\in \mathbb Z^2\setminus (0,0)$
$$
d^*(A+B\cap A+B-(n,m)\cap A+B-2(n,m) \cap A+B-3(n,m) \cap A+B-4(n,m))\leq \frac{1}{2}d^*(B)^5.
$$
If so, can the exponent $5$ be improved?
\end{question}

The matter of optimal recurrence in $\mathbb Z^d$ is still under investigation; for configurations of 3 points in an arbitrary set of positive density, the best known results are due to Q. Chu (\cite{Chu}).

\subsection{Characteristic factors}

If $p_i:\mathbb Z\to \mathbb Z,i=1,\dots, k-1$ are polynomials and $\mathbf X=(X,\mathscr X,\mu,T)$ is a M.~P.~S.~ with $f\in L^2(\mu),$ consider the correlation sequence
$$
I_{\vec{p}}(f;n)= \int f\cdot f\circ T^{p_1(n)} \cdot f\circ T^{p_2(n)}\cdot \cdots \cdot f\circ T^{p_{k-1}(n)}\, d\mu.
$$
For a given ergodic $\mathbf X,$ there is a factor $\mathbf Z=(Z,\mathscr Z,\nu,S)$ of $\mathbf X$ such that the sequences $I_{\vec{p}}(f;n)$ are largely unaffected by replacing $f$ with $\mathbb E(f|\mathscr Z).$  Furthermore, $\mathbf Z$ is an inverse limit of nilsystems.

\begin{proposition}[\cite{BHK}, Corollary 4.6, cf.~\cite{Lei?}]\label{nilfactors}  Let $\mathbf X$ be an ergodic M.~P.~S.~ and let $p_1,\dots,p_{k-1}:\mathbb Z\to \mathbb Z$ be polynomials.  Then there is a factor $\mathbf Z$ of $\mathbf X$ such that
\begin{align}\label{tofactor}
d^*\{n:|I_{\vec{p}}(f;n)-I_{\vec{p}}(\mathbb E(f|\mathscr Z);n)|>\varepsilon\}=0
\end{align}
for all $\varepsilon>0.$  Furthermore, $\mathbf Z$ is an inverse limit of nilsystems.
\end{proposition}

The combination of the above theorem with the following lemma is essentially the proof of Theorem \ref{main}.

\begin{lemma}\label{cesaro} Let $\mathbf Z$ be an inverse limit of ergodic nilsystems, with Kronecker factor $\mathbf Z_1,$ and let $A\subset \mathbb Z$ be infinite.  If $f\in L^2(\mu)$ is orthogonal to $L^2(\mathscr Z_1),$ then there is a sequence $\{a_n\}_{n\in \mathbb N}$ of elements of $A$ such that
$$
\lim_{N\to \infty} \frac{1}{N}\sum_{n=1}^N f\circ T^{-a_n}=0
$$
in $L^2(\mu).$
\end{lemma}

\textit{Proof.} We use the following fact: if $(x_n)_{n\in \mathbb N}$ is a bounded sequence of vectors in a Hilbert space such that $x_n\to 0$ weakly, then there is a subsequence $(x_n')_{n\in \mathbb N}$ such that $\lim_{N\to \infty} \bigl\|\frac{1}{N}\sum_{n=1}^N x_n\bigr\|=0.$

By Theorem \ref{lebesguespectrum}, $\lim_{n\to \infty} \int f\circ T^n \cdot g\,d\mu=0$ for all $g\in L^2(\mu).$  In particular, we can choose a sequence of elements $\{a_n'\}_{n=1}^\infty$ of $A$ with $\lim_{n\to \infty} \int f\circ T^{-a_n'}\cdot g\, d\mu=0$ for all $g\in L^2(\mu),$ and a subsequence $\{a_n\}_{n\in \mathbb N}$ of $\{a_n'\}_{n\in \mathbb N}$ such that $\lim_{N\to \infty} \frac{1}{N}\sum_{n=1}^N f\circ T^{-a_n}=0$ in norm.  \hfill $\square$

Lemma \ref{cesaro} leads immediately to a proof of Theorem \ref{main'}, specialized to inverse limits of nilsystems.

\begin{theorem}\label{specialcase}  Let $\mathbf Z$ be an inverse limit of ergodic nilsystems, $f:Z\to [0,1].$ Let $p_1,\dots, p_{k-1}:\mathbb Z\to \mathbb Z$ be polynomials such that $\mathbb N$ is a set of $\vec{p}$-recurrence, and let $A\subset \mathbb Z$ be infinite. Then for all $\varepsilon>0,$ there exists a finite $A'\subset A$ with $f_{A'}:=\frac{1}{|A'|}\sum_{a\in A'} f\circ T^{-a}$ such that
$$
\bigl\{n:I_{\vec{p}}(f_{A'};n)> \Bigl(\int f\, d\mu \Bigr)^k-\varepsilon\bigr\}
$$
is $\text{Nil}_r$-Bohr, where $r=\deg{\vec{p}}.$
\end{theorem}

\textit{Proof.}   First note that we can assume that $\mathbf Z$ is actually a nilsystem, since we can approximate $f$ in $L^2(\mu)$ by functions which are measurable with respect to some nilsystem factor of $\mathbf Z.$  Applying Cauchy-Schwarz and the triangle inequality, we see that $|I_{\vec{p}}(f)-I_{\vec{p}}(g)|<\delta$ whenever $\|f-g\|_{L^2(\mu)}<\delta/2k.$

Let $\mathbf Z_1=(Z_1,\mathscr Z_1,m_1,R_{\alpha})$ be the Kronecker factor of $\mathbf Z.$  Write $f=g+h,$ where $\mathbb E(g|\mathscr Z_1)=0,$ and $h$ is $\mathscr Z_1$-measurable.    By Lemma \ref{cesaro}, there exists a sequence $\{a_n\}_{n\in \mathbb N}$ of elements of $A$ satisfying $\lim_{N\to \infty} \frac{1}{N}\sum_{n=1}^N g\circ T^{-a_n}=0$ in $L^2(\mu).$  Let $A'=\{a_n:1\leq n\leq N\}$ for some $N$ satisfying $\|\frac{1}{N}\sum_{n=1}^N g\circ T^{-a_n}=0\|< \varepsilon/2k.$  Then $f_{A'}=g_{A'}+h_{A'},$ where $\|g_{A'}\|<\varepsilon/2k.$  Thus
\begin{align}\label{close}
|I_{\vec{p}}(f_{A'};n)-I_{\vec{p}}(h_{A'};n)|\leq \varepsilon/2
\end{align} for all $n.$  Since $h_{A'}$ is $\mathscr Z_1$-measurable, we can interpret $I_{\vec{p}}(h_{A'};n)$ as $I_{\vec{p}}(\tilde h_{A'};n),$ where $\tilde h_{A'}: Z_1\to [0,1]$ is measurable, and $\tilde h_{A'}(T^{p_i(n)}z)=h_{A'}(z+p_i(n)\alpha).$  Thus
\begin{align}\label{Kintegral}
I_{\vec{p}}(\tilde h_{A'};n) = \int \tilde{h}_{A'}(z)\cdot \tilde{h}_{A'}(z+p_1(n)\alpha)\cdot \cdots \tilde{h}_{A'}(z+p_{k-1}(n)\alpha)\, d\mu.
\end{align}
Since $\vec{p}$ is a recurrent set of polynomials, we know that, for every neighborhood $U$ of  $0\in Z_1,$ the set
$$
\{n:p_{i}(n)\alpha \in U \text{ for } i=1,\dots, k\}
$$
Is $\text{Nil}_r$-Bohr.  For a sufficiently small neighborhood $U, p_i(n)\alpha\in U$ implies that the integral in (\ref{Kintegral}) is at least $\int \bigl( h_{A'}\bigr)^k \, d\mu-\varepsilon/2.$  Since $\int \bigl( h_{A'}\bigr)^k \, d\mu\geq \bigl(\int h_A \, d\mu\bigr)^k = \bigl(\int f\, d\mu\bigr)^k,$ equations (\ref{close}) and (\ref{Kintegral}) imply the desired conclusion.  \hfill $\square$

\medskip

\textit{Proof of Theorem \ref{main'}.}  Let $\mathbf X$ be an arbitrary ergodic system, let $\vec{p}$ be a polynomial vector as in the hypothesis, and let $\mathbf Z$ be the factor of $\mathbf X$ guaranteed to exist by Proposition \ref{nilfactors}.  Let $f:X\to [0,1],$ let $h=\mathbb E(f|\mathscr Z),$ and choose $A'\subset A$ so that $R:=\{n:I_{\vec{p}}(h_{A'}; n)\geq \bigl(\int f\, d\mu\bigr)^k-\varepsilon/2\}$ is $\text{Nil}_r$-Bohr, as Theorem \ref{specialcase} allows.  By (\ref{tofactor}), $|I_{\vec{p}}(h_{A'};n)-I_{\vec{p}}(f_{A'};n)|<\varepsilon$ on a set $R'=R\setminus E,$ where $d^*(E)=0.$  Hence $I_{\vec{p}}(f_{A'},n)> \bigl(\int f\, d\mu\bigr)^k$ for $n\in R'.$

By hypothesis, the given set $S$ is a frequent set of $\vec{p}$-recurrence, so $d^*(S\cap R)>0,$ and $d^*(S\cap R')>0,$ as was to be shown.  Since $R'$ is almost $\text{Nil}_r$-Bohr, we are done.  \hfill $\square$

\bigskip

\section{Proof of theorem \ref{near}}  Given an M.~P.~S.~ $(X,\mathscr X,\mu,T),$ we consider the correlation sequences $\int f\circ T^n\cdot g.$  By Herglotz's theorem, one can write
$$
\int f\circ T^n \cdot g\, d\mu= \int e^{2\pi in\theta}\, d\sigma(\theta)
$$ where $\sigma$ is a complex measure on $\mathbb T.$ The eigenvalues of $T$ correspond to the atoms of $\sigma,$ so that $\sigma$ is atomless whenever $f$ is orthogonal to all eigenfunctions of $T.$

\medskip

The following variation on Wiener's lemma establishes estimates on $\hat{\sigma}(n): \int e^{2\pi i n\theta}\, d\sigma(\theta).$

\medskip

\begin{lemma}\label{Wiener}
Let $\sigma$ be an atomless complex measure on $[0,1],$ and suppose that $\mathcal S=(S_j)_{j\in \mathbb N}$ is a nearly equidistributed sequence of finite subsets of $\mathbb Z.$  Then for all $\varepsilon>0,$
$$
d_{\mathcal S}\{n: |\hat{\sigma}(n)|>\varepsilon\}=0.
$$
\end{lemma}
The proof is essentially the same as the proof of Wiener's lemma when $S_j=\{1,\dots, j\}.$

\medskip

\textit{Proof of Theorem \ref{near}.} Let $(X,\mathscr X,\mu, T)$ be an ergodic M.~P.~S. with Kronecker factor $\pi: \mathbf X\to \mathbf Z=(Z,\mathscr Z,m,R_\alpha),$ and let $D\in \mathscr X.$   Let $(S_j)_{j\in \mathbb N}$ be a nearly equidistributed sequence, and let $A\subset \mathbb Z$ with $d_{\mathcal S}(A)>0.$  We will show that $D_A:=\bigcup_{a\in A} T^{-a}D$ contains a set of the form $\pi^{-1}(K),$ where $K\subset Z$ and $m(K)\geq \mu(D).$  Thus,
$$
\mu(D_A\cap T^{-n_1}D_A\cap \dots \cap T^{-n_k}D_A)\geq m(K\cap K-n_1\alpha \cap K-n_2\alpha \cap \dots \cap K-n_{k-1}\alpha)
$$
The conclusion then follows from the fact that the map $Z^{k-1}\to \mathbb R$ given by $(z_1,\dots,z_{k-1})\mapsto m(K\cap K-z_1\cap K-z_2\cap \dots \cap K-z_{k-1})$ is continuous in $(z_1,\dots, z_k).$

\medskip

Write $f=1_D,$ and write $f= g+h,$ where $g=\mathbb E(f|\mathscr Z), h\perp L^2(\mathscr Z).$    We will estimate $1_{D_A}$ from below, noting that for all finite $A'\subset A, f_{A'}:=\frac{1}{|A'|}\sum_{a\in A'} f\circ T^{-a}\leq 1_{D_A}$ (pointwise).  In particular, if $f_0$ is a weak limit of the averages $f_{A'},$ then $f_0\leq 1_{D_A}.$  We will show that some weak limit $f_0$ satisfies $f_0(x)>0$ for some $x\in \pi^{-1}(K),$ where $K\subset Z$ with $m(K)\geq \mu(D).$

\medskip

Now $h$ is orthogonal to the eigenfunctions of $T,$ so for a given $\phi\in L^2(\mu)$ we can write
$$
\int h\circ T^n \cdot \phi, d\mu = \int e^{2\pi i n\theta}\, d\sigma
$$
where $\sigma$ is an atomless complex measure on $[0,1].$  Choose a subsequence $(S_j')_{j\in \mathbb N}$ from $\mathcal S$ so that $\lim_{j\to \infty} \frac{|A\cap S_j|}{|S_j|}=d_{\mathcal S}(A),$ and write $A_j$ for $A\cap S_j'.$  By Lemma \ref{Wiener}, we have that for all $\varepsilon>0,$
$$
d_{\mathcal S}\Bigl\{n:\bigl|\int h\circ T^n \cdot \phi\, d\mu\bigr|>\varepsilon\Bigr\}=0,
$$
so
$$
\lim_{j\to\infty} \frac{1}{|A_j|}\sum_{a\in A_j} \int h\circ T^{-a}\cdot \phi\, d\mu=0.
$$
By a diagonalization procedure, we may now pick a sequence $a_n$ of elements of $A$ so that
\begin{enumerate}
\item[(i)] $\lim_{N\to \infty} \frac{1}{N}\sum_{n=1}^N h\circ T^{-a_n}=0$ (weakly in $L^2(\mu)$).
\item[(ii)] $f_0:=\lim_{N\to \infty} \frac{1}{N}\sum_{n=1}^N f\circ T^{-a_n}$ exists in $L^2(\mu).$
\end{enumerate}
Since $f$ is nonnegative and bounded above by $1,$ so is $f_0,$ and $\int f_0, d\mu\geq \int f\, d\mu =\mu(D).$  Hence $f_0$ is positive on a set of measure at least $\mu(D).$  By (i) and the splitting $f=g+h,$ the limit $f_0$ is measurable with respect to $\pi^{-1}(\mathscr Z).$  We may now take $K=\{x:f_0(x)>0\}.$  \hfill $\square$

\bigskip

\section{Appendix.}\label{appendix}

Let $G$ be a nilpotent Lie group, and let $\Gamma$ be a discrete subgroup of $G$ such that $X:=G/\Gamma$ is compact.  Let $\tau\in G,$ and let $T:X\to X$ be the translation $Tx\Gamma=\tau x\Gamma,$ which preserves the natural projection of Haar measure $\mu$ on $X.$  Suppose that $T$ is ergodic.  We want to show that when $f$ is orthogonal to the eigenfunctions of $T,$ then $(X,\mathscr X,\mu,T),$ then $\lim_{n\to \infty} \int f\circ T^n\cdot g\, d\mu=0$ for all $g\in L^2(\mu).$  We derive this as a corollary of the following result.

\begin{theorem}[cf.~\cite{Parry}, Theorem 3]\label{Parrythm}  Let $H$ be a connected, simply connected nilpotent Lie group, with $\Gamma\subset H$ a discrete subgroup and $X:=H/\Gamma$ compact with Haar measure $\mu$.  Let $A:H\to H$ be a unipotent automorphism and $a\in H.$  Suppose that the map $S:X\to X$ given by $Sh\Gamma=aA(h)\Gamma$ is ergodic.  Then the operator given by $f\mapsto f\circ S$ has countable Lebesgue spectrum in the orothocomplement of the space spanned by the eigenfunctions of $S.$  In particular if $f$ is orthogonal to the eigenfunctions of $S,$ then $\int f\circ S^n\cdot g\, d\mu\to 0$ as $n\to \infty.$
\end{theorem}
Here ``unipotent" means the map $B(h):=A(h)h^{-1}$ satisfies $B^n(h)=e_G,$ for some $n$ and all $h.$

\medskip

Assuming $G$ is as above, Let $G^o$ be the connected component of $G.$  Passing to the universal cover of $G^o,$ we may assume $G^o$ is simply connected.  When $X$ is connected, we have $G=G^o\Gamma,$ so $\tau$ can be written as $t\gamma,$ where $t\in G^o$ and $\gamma\in \Gamma.$  Also $G/\Gamma \equiv G^o/(\Gamma\cap G^o).$  Writing
$$
Tx\Gamma= t\gamma x \Gamma = t \gamma x \gamma^{-1} \Gamma,
$$  we can consider $T$ as an action on $X':=G^o/(\Gamma \cap G^o).$  This action has the form $x\Gamma\mapsto tAx\Gamma,$ where $A$ is a unipotent automorphism of $G^o$ (note that all inner automorphisms of a nilpotent Lie group are unipotent).  We can apply Theorem \ref{Parrythm} to obtain our conclusion.

\medskip

In general, $X$ may not be connected, but $X$ will have finitely many connected components $X_1,\dots, X_r$ which are permuted by $T.$  In this case, the map $T^{r!}:X_1\to X_1$ will be ergodic, so we may apply the above argument to $T^{r!}.$  A simple argument then shows that the desired conclusion holds for $T$ as well.

\end{document}